\documentclass[11pt,a4paper]{article}

\usepackage{header}
\usepackage{abbrev}

\begin{document}

% Titolo
\title{A single shooting method with approximate Fr\'{e}chet derivative for computing geodesics on the Stiefel manifold}

% Authors: full names plus addresses.
\author{Marco Sutti\thanks{Mathematics Division, National Center for Theoretical Sciences, National Taiwan University, No. 1, Sec. 4, Roosevelt Road, Taipei 10617, Taiwan (\email{msutti@ncts.ntu.edu.tw}).}\hspace{2mm}\orcidlink{0000-0002-8410-1372}}

\date{\today}

\maketitle

\begin{abstract}
This paper shows how to use the shooting method, a classical numerical algorithm for solving boundary value problems, to compute the Riemannian distance on the Stiefel manifold $\Stnp$, the set of $ n \times p $ matrices with orthonormal columns. The proposed method is a shooting method in the sense of the classical shooting methods for solving boundary value problems; see, e.g., Stoer and Bulirsch, 1991.
The main feature is that we provide an approximate formula for the Fr\'{e}chet derivative of the geodesic involved in our shooting method.
Numerical experiments demonstrate the algorithms' accuracy and performance. Comparisons with existing state-of-the-art algorithms for solving the same problem show that our method is competitive and even beats several algorithms in many cases.

%Include keywords and mathematical subject classification numbers as needed.
\bigskip
\textbf{Key words.} Stiefel manifold, shooting methods, endpoint geodesic problem, Riemannian distance, Newton's method, Fr\'{e}chet derivative

\medskip
\textbf{AMS subject classifications.} 65L10, 65F45, 65F60, 65L05, 53C22, 58C15

\end{abstract}

\section{Introduction}

The object of study in this paper is the compact Stiefel manifold, i.e., the set of $ n \times p $ matrices with orthonormal columns
\[
   \Stnp = \left\lbrace X \in \R^{n \times p}: \ X\tr\! X = I_p \right\rbrace.
\]
There are applications in several areas of mathematics and engineering that deal with data that belong to $\Stnp$. Domains of applications include numerical optimization, imaging, and signal processing. Some applications, like finding the Riemannian center of mass, require evaluating the geodesic distance between two arbitrary points on $\Stnp$. Since no explicit formula is known for computing the distance on $\Stnp$, one has to resort to numerical methods. 

In this paper, we are concerned with computing the Riemannian distance between two given points on the Stiefel manifold. As we shall see, the distance between two points on a manifold is related to the concept of minimizing geodesic\footnote{Geodesics are generally defined as critical points of the length functional, and as such, they may or may not be minima. A minimizing geodesic is one that minimizes the length functional. We introduce the notion of geodesics in Section~\ref{sec:geodesic_exp_log}.}.
The problem can be briefly formulated as follows. Given two points $X$, $Y$ on $\Stnp$ that are sufficiently close to each other, finding the distance between them is equivalent to finding the tangent vector $\xi^{\ast} \in \mathrm{T}_{X}\Stnp$ with the shortest possible length such that \protect{\cite{Lee:2018,boumal2023intromanifolds}}
\[
\Exp_{X}(\xi^{\ast}) = Y,
\]
where $\Exp_{X}$ denotes the Riemannian exponential mapping at $X$. 
The solution to this problem is equivalent to the Riemannian logarithm of $Y$ with base point $X$
\[
 \xi^{\ast} = \Log_{X}(Y).
\]
The sought distance between $X$ and $Y$ is then given by the norm of $ \xi^{\ast} $.

\begin{figure}[htbp]
\centering
\includegraphics[width=0.50\columnwidth]{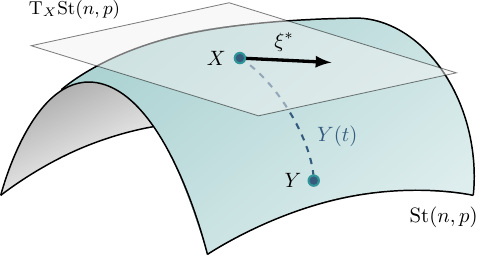}
\caption{Illustration of the problem statement.}
\label{fig:single_shooting}
\end{figure}

Figure~\ref{fig:single_shooting} provides an artistic illustration of the problem. The latter will be stated in more detail in Section~\ref{sec:problem_statement}.

%=========================================
% FIGURE: PROBLEM STATEMENT
%=========================================
%=========================================

It is interesting to note that, for some manifolds, explicit formulas exist for computing the Riemannian distance. This is the case of the Grassmann manifold $ \Grass(n,p) $, which is the set of $ p $-dimensional vector subspaces of $ \R^{n} $. For instance, let $ \cX $ and $ \cY $ be two subspaces belonging to $ \Grass(n,p) $, then the distance between $ \cX $ and $ \cY $ is
\[
   d(\cX,\cY) = \sqrt{\theta_{1}^{2} + \cdots + \theta^{2}_{p}},
\]
where $ \theta_{i} $, $ i = 1, \ldots, p $, are the principal angles between $ \cX $ and $ \cY $; see, e.g., \protect{\cite[Theorem~8]{Wong:1967}} and \protect{\cite[Section~3.8]{AMS:2004}}. The unit sphere $ \cS^{n-1} $ embedded in $\R^{n}$ also has explicit formulas for computing the Riemannian distance. In contrast, no such closed-form solution is known for the Stiefel manifold. This motivates us to consider numerical methods. In general, the problem of finding the distance given two points on a Riemannian manifold is related to the Riemannian logarithm function (more details later in Section~\ref{sec:geodesic_exp_log}). Several authors have already tackled the problem of computing the Riemannian logarithm on the Stiefel manifold. These contributions are detailed in Section~\ref{sec:other_approaches}.

\subsection{Contributions}

In this work, we use the shooting method, which is a classical numerical algorithm for solving boundary value problems, to compute the distance on the Stiefel manifold $\Stnp$. These methods are not new (thorough coverage of the shooting methods is given, e.g., in~\protect{\cite{Stoer:1991}}), but their application to computing the Riemannian distance on the Stiefel manifold is relatively new. The method of Bryner, although also named ``shooting'' in~\protect{\cite{Bryner:2017}}, cannot be regarded as a classical shooting method since it makes use of Riemannian geometry concepts (like the parallel transport) that do not fit in the classical framework of~\protect{\cite{Stoer:1991}}.

In an earlier version of this work~\protect{\cite{Sutti:2020b,Sutti:2023}}, we used the vectorization operator and Kronecker products to work out explicit expression for the Jacobian matrices involved in the shooting method. The reason why this was done was to carry out some preliminary analysis on the explicit expressions of the Jacobian matrices. The drawback in the actual implementation was the excessive computational cost given the dimensions of the operators involved, even for small values of $n$ and $p$.

Here, in contrast, we work directly with matrices, and we use an approximate form of the Fr\'{e}chet derivative given by a truncated Fr\'{e}chet derivative of the matrix exponential. Hence, there is no need for finite difference approximations. In particular, the main contributions of this paper are as follows.

\begin{enumerate}[label=(\roman*)]
\item We provide a single shooting method for computing geodesics on the Stiefel manifold using the canonical metric as a classical numerical algorithm for solving boundary value problems.
\item We introduce a truncated Fr\'{e}chet derivative that leads to a linear matrix equation that can be efficiently solved to find the algorithmic update.
\item We perform extensive numerical experiments to demonstrate the algorithms in terms of performance and accuracy. In particular, comparisons with existing state-of-the-art algorithms for solving the same problem show that our method is competitive and even beats several algorithms in many cases.
\end{enumerate}

\subsection{Applications and motivation}

Many scientific and engineering works have used the Stiefel manifold in their applications. To provide some motivation for the present work, this section summarizes a few applications that explicitly compute the geodesic distance.

In affine invariant shape analysis, Younes et al.~\protect{\cite{Younes:2008}} studied a specific metric on plane curves that has the property of being isometric to classical manifolds (like the sphere, complex projective plane, Stiefel and Grassmann manifolds) modulo change of parametrization. Moreover, they provided experimental results that explicitly compute minimizing geodesics between two closed curves.

In the context of shape analysis of closed curves, Srivastava and Klassen \protect{\cite{Srivastava:2016}} studied the space of functions representing unit-length, planar, closed curves, which can be shown to be a Stiefel manifold.
Ring and Wirth \protect{\cite[Section~4.2]{Ring:2012}} provided an application for image segmentation on the Stiefel manifold using a Riemannian variant of the classical BFGS algorithm. This is compared to the work of Sundaramoorthi et al.~\protect{\cite{Sundaramoorthi:2011}}, where the authors used geodesic retractions based on the matrix exponential.
The more general reference by Kendall et al.~\protect{\cite[Chapter~6]{Kendall:1999}} also contains a discussion on the Stiefel manifold and shape spaces. 
Bryner \protect{\cite{Bryner:2017}} also proposed some numerical applications on the pre-shape space.
 
\c {C}eting\"ul and Vidal~\protect{\cite{Cetingul:2009}} investigated the intrinsic mean shift algorithm for clustering on Stiefel and Grassmann manifolds. Turaga et al.~\protect{\cite{Turaga:2008,Turaga:2011}} investigated applications of the Stiefel manifold in computer vision and pattern recognition to develop accurate inference algorithms. Vision applications such as activity recognition, video-based face recognition, shape classification, and unsupervised clustering were targeted. In particular, step 3 of Algorithm 1 in~\protect{\cite{Turaga:2011}} computes the inverse exponential map, but it was unclear how this was achieved.

The low-rank representation (LRR) is a widely used technique in computer vision and pattern recognition for data clustering models. Yin et al.~\protect{\cite{Yin:2015}} extended the LRR from Euclidean space to the manifold-valued data on the Stiefel manifold by incorporating the intrinsic geometry of the manifold. They acknowledged that, in general, it is pretty hard to compute the log mapping for the Stiefel manifold. Consequently,  they used the retraction map (a first-order approximation to the exponential mapping~\protect{\cite{Absil:2012}}) instead of the exponential map because of its reduced computational cost.

More recently, Li and Ma~\protect{\cite{LiMa:2022}} proposed a generalization of the federated learning framework to Riemannian manifolds. In particular, they consider the kPCA problem on the Stiefel manifold. Even though they initially discuss the Riemannian logarithm mapping, they finally adopt a retraction in the numerical implementations, similarly to what was done by \protect{\cite{Yin:2015}}.

\subsection{Related works and other approaches} \label{sec:other_approaches}

Shooting methods are not the only option to solve the endpoint geodesic problem; many other numerical algorithms have been proposed. As a plethora of methods now come out every year, this review is not meant to be exhaustive.

The leapfrog algorithm by Noakes~\protect{\cite{Noakes:1998}} is based on partitioning the original problem into smaller subproblems. This method has global convergence properties, but it slows down for an increasing number of subproblems or when the solution is approached \protect{\cite[Section~1]{Kaya:2008}}. Moreover, Noakes realized that his leapfrog algorithm was in some way imitating the Gauss--Seidel method \protect{\cite[Section~1]{Noakes:1998}}. This connection has been explored by Sutti and Vandereycken~\protect{\cite{Sutti_V:2023}}.

Bryner~\protect{\cite{Bryner:2017}} proposed two schemes, named ``shooting method'' and ``path-straightening'', to compute endpoint geodesics on the Stiefel manifold by considering them as an embedded submanifold of the Euclidean space. From the matrix algebra perspective, Rentmeesters~\protect{\cite{Rentmeesters:2013}} and Zimmermann~\protect{\cite{Zimmermann:2017,Zimmermann:2019}} derived algorithms for evaluating the Riemannian logarithm map on the Stiefel manifold with respect to the canonical metric, which is locally convergent and depends upon the definition of the matrix logarithm function. Recently, Zimmermann and H\"{u}per~\protect{\cite{Zimmermann:2022}} provided a unified method to deal with the endpoint geodesic problem on the Stiefel manifold with respect to a family of metrics.

Noakes and Zhang~\protect{\cite{Noakes:2022}} proposed an alternative algorithm to find geodesics joining two given points. Like leapfrog, this method exploits the shooting method to compute geodesics joining junction points.

The methods of Nguyen~\protect{\cite{Nguyen:2022}} are based on classical (Euclidean) optimization algorithms for minimizing an objective function with a reduction of the computational cost thanks to the formulation of the gradients using only the Fr\'{e}chet derivatives.

\subsection{Notation}\label{sec:notation}

Here, we list the notations and symbols adopted in this paper in order of appearance. Symbols only used in one section are typically omitted from this list.

%=====================================================
% TABLE: LIST OF SYMBOLS
%=====================================================
\begin{table}[htbp]
   %\begin{center}
      \begin{tabular}{ll} 
          $ \Stnp $                          &  Stiefel manifold of orthonormal $n$-by-$p$ matrices \\
          $ X $, $ Y $, $ Y_{0} $, $ Y_{1} $ &  Elements of $ \Stnp $ \\
          $ I_{p} $                          &  The identity matrix of size $p$-by-$p$ \\
          $ \mathrm{T}_{X}\Stnp $            &  Tangent space at $ X $ to the Stiefel manifold $ \Stnp $ \\
          $ \xi^{\ast} $                     &  A tangent vector that we want to recover \\
          $ \Exp_{X} $                       &  Riemannian exponential map at $X$ \\
          $ \Log_{X} $                       &  Riemannian logarithm map at $X$ \\
          $ O_{p} $                          &  The null matrix of size $p$-by-$p$ \\
          $ I_{n,p} $                        &  The matrix $ \begin{bmatrix} I_{p} \\ O_{(n-p)\times p} \end{bmatrix} $ \\
          $ \cS_{\mathrm{sym}}(p) $          &  Space of $ p $-by-$ p $ symmetric matrices \\
          $ \cS^{n-1} $                      &  The unit sphere embedded in $\R^{n}$ \\
          $ \On $                            &  The orthogonal group of $n$-by-$n$ orthogonal matrices \\
          $ X_{\perp} $                      &  An orthonormal matrix whose columns span \\
                                             &  the orthogonal complement of $ \mathrm{span}(X) $ \\
          $ \cS_{\mathrm{skew}}(p) $         &  Space of $ p $-by-$ p $ skew-symmetric matrices \\ 
          $ \Omega $                         &  An element of $ \cS_{\mathrm{skew}}(p) $ \\
          $ K $                              &  A matrix in $ \R^{(n-p) \times p} $ \\
          $ \cM $                            &  Generic manifold \\
          $ \mathrm{T}_{x}\cM $              &  Tangent space at $ x $ to the manifold $ \cM $ \\
          $ \langle \cdot, \cdot \rangle_{x} $   &  Inner product on the tangent space $ \mathrm{T}_{x}\cM $ \\
          $ g $                              &  Riemannian metric \\
          $ \gamma(t) $                      &  Parametrized curve on the manifold $\cM$ \\
          $ d $                              &  Riemannian distance function \\
          $ d(x,y) $                         &  Riemannian distance between two points $x$ and $y$ \\
          $ \P_{X} $                         &  The projector onto the tangent space $ \mathrm{T}_{X}\Stnp $ \\
          $ \inj_{X}(\cM) $                  &  Injectivity radius of $\cM$ at $X$ \\
          $ \inj(\cM) $                      &  Global injectivity radius of $\cM$ \\
          $ \| \cdot \|_{\F} $               &  Frobenius norm \\
          $ \| \cdot \|_{\mathrm{c}} $       &  Canonical norm \\
          $ \expm $              &  The matrix exponential \\
          $ A $   & The matrix $ \begin{bmatrix}   \Omega   &   -K\tr \\   K    &    O_{n-p}  \end{bmatrix} $ \\
          $ Z_{1}(t) $ or $ Y(t) $           &  A geodesic on $\Stnp$ \\
          $ Z_{2}(t) $ or $ \dt{Y}(t) $      &  The derivative of a geodesic \\
          $ F $                              &  The nonlinear function $ Z_{1}(1,\xi) - Y_{1} $ \\
          $ F^{(k)} $                        &  The nonlinear function $F$ evaluated at iteration $ k $ \\
          $ \delta \xi^{(k)} $               &  The residual at iteration $ k $ in the single shooting method \\
          $ \D \! f(A)[E] $                  &  Fr\'{e}chet derivative of a matrix function $f$ at $A$ in the direction $E$
      \end{tabular}
   %\end{center}
\end{table}
%=========================================

\subsection{Outline of the paper}

The remaining part of this paper is organized as follows.
Section~\ref{sec:stiefel_manifold} introduces the geometry of the Stiefel manifold. Readers who are familiar with Riemannian geometry, particularly the geometry of the Stiefel manifold, might want to skip this section.
Section~\ref{sec:problem_statement} presents the problem statement, which is the focus of this work. Section~\ref{sec:SSAF} describes our proposed algorithm: a single shooting method with an approximate Fr\'{e}chet derivative. Numerical experiments and comparisons with other methods are presented in Section~\ref{sec:numerical_experiments}. Finally, we conclude the paper by summarizing the contributions and providing future research outlooks in Section~\ref{sec:conclusions}.

\section{Geometry of the Stiefel manifold}\label{sec:stiefel_manifold}

This section introduces the geometry of the Stiefel manifold. Here, we only give the necessary background to understand the remaining part of this paper. For additional details, we refer the reader to the reference works of~\protect{\cite{Edelman:1998,AMS:2008,boumal2023intromanifolds}}.

The set of all $ n \times p $ orthonormal matrices, i.e.,
\[
   \Stnp = \lbrace X \in \Rnp \colon X\tr \! X = I_{p} \rbrace,
\]
endowed with its submanifold structure is called an orthogonal or compact Stiefel manifold. It is a subset of $ \Rnp $, and it can be proven that it has the structure of an embedded submanifold of $ \Rnp $ \protect{\cite[Section~3.3.2]{AMS:2008}}. 
The Stiefel manifold $ \Stnp $ may also degenerate to some special cases. For $ p = 1 $, it reduces to the unit sphere $ \cS^{n-1} $ in $ \R^{n} $, while for $ p = n $, it becomes the orthogonal group $ \On $, whose dimension is $ \tfrac{1}{2}n(n-1) $.

\subsection{Tangent spaces and projectors} \label{sec:tgspaces_stiefel}

The tangent space to a manifold at a given point can be seen as a local vector space approximation to the manifold at that point. In practice, the tangent space is used to perform the operations of vector addition and scalar multiplication, which would otherwise be impossible to perform on the manifold without leaving it due to the manifold's curvature. Endowed with a Euclidean inner product, this vector space becomes a Euclidean space where we also have a notion of lengths. Here, we will directly focus on the tangent space to the Stiefel manifold. For a more precise definition of a tangent space in the general case, we refer the reader to \protect{\cite{AMS:2008}}.

The tangent space to the Stiefel manifold at a point $X$ can be characterized by \protect{\cite[Example 3.5.2]{AMS:2008}}
\begin{equation}\label{eq:tg_space_stiefel}
   \mathrm{T}_{X}\Stnp = \lbrace X \Omega + X_{\perp} K \colon  \Omega = -\Omega\tr, \ K \in \R^{(n-p)\times p} \rbrace,
\end{equation}
$ \Omega $ being a $p$-by-$p$ skew-symmetric matrix, $ \Omega \in \cS_{\mathrm{skew}}(p) $, $ X_{\perp} $ being an orthonormal matrix whose columns span the orthogonal complement of $ \mathrm{span}(X) $, and $ K \in \R^{(n-p) \times p} $, with no restriction on $ K $. With this characterization in mind, and with the fact that $ \dim\!\big(\Stnp\big) = \dim\!\big(\mathrm{T}_{X}\Stnp\big) $, it is straightforward to calculate the dimension of the Stiefel manifold as
\[
   \dim(\Stnp) = \dim(\cS_{\mathrm{skew}}(p)) + \dim(\R^{(n-p) \times p}) = \tfrac{1}{2}p(p-1) + (n-p) p = np - \tfrac{1}{2}p(p+1).
\]
The projection onto the tangent space $ \mathrm{T}_{X}\Stnp $ is
\begin{equation}\label{eq:proj_tg_space}
   \P_{X}\xi = X \mathrm{skew}(X\tr \xi) + (I-XX\tr)\,\xi,
\end{equation}

\subsection{Geodesics, exponential mapping and logarithm mapping} \label{sec:geodesic_exp_log}

\begin{figure}[htbp]
\centering
\includegraphics[width=0.50\columnwidth]{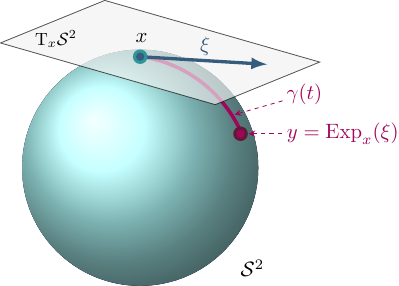}
\caption{The Riemannian exponential map on the sphere.}
\label{fig:exponential_map}
\end{figure}

Geodesics are defined in general as curves with zero covariant acceleration. They allow us to introduce the \emph{Riemannian exponential} $ \Exp_{x}\colon \mathrm{T}_{x}\cM \to \cM $ that maps a tangent vector $\xi = \dt{\gamma}(0) \in \mathrm{T}_{x}\cM $ to the geodesic endpoint $\gamma(1) = y $: $ \Exp_{x}(\xi) = y$. Figure~\ref{fig:exponential_map} illustrates these concepts for the unit sphere $ \cS^{2} $, which is also a special case of Stiefel manifold $\Stnp$, with $n=3$ and $p=1$.

%=========================================
% FIGURE: EXPONENTIAL MAP
%=========================================
%=========================================

To define a distance on a given manifold $\cM$, we need a notion of length that applies to tangent vectors. To this aim, we endow the tangent space $ \mathrm{T}_{x}\cM $ with an \emph{inner product} $ \langle \cdot, \cdot \rangle_{x} $.
The inner product $ \langle \cdot, \cdot \rangle_{x} $ induces a norm $ \| \xi_{x} \|_{x} = \sqrt{\langle \xi_{x}, \xi_{x} \rangle_{x}} $ on $ \mathrm{T}_{x}\cM $.
A manifold $ \cM $ endowed with a smoothly-varying inner product (called \emph{Riemannian metric} $ g $) is called \emph{Riemannian manifold}.

The \emph{length of a curve} $ \gamma \colon [a,b] \to \cM $ on a Riemannian manifold $ ( \cM, g ) $ is
\[
   L(\gamma) = \int_{a}^{b} \sqrt{g(\dt{\gamma}(t),\dt{\gamma}(t))} \, \mathrm{d}t.
\]

The \emph{Riemannian distance} is defined as the shortest path between two points $ x $ and $ y $
\[
   d \colon \cM \times \cM \to \R \colon d(x,y) = \inf_{\gamma \in \Gamma} L(\gamma),
\]
where $ \Gamma $ denotes the set of all curves $ \gamma $ in $ \cM $ joining points $ x $ and $ y $.

Generally speaking, different choices of Riemannian metric are possible. In this paper, we consider the non-Euclidean \emph{canonical metric} inherited by $\Stnp$ from its definition as a quotient space of the orthogonal group \protect{\cite[(2.39)]{Edelman:1998}}. Given $ X \in \Stnp $ and $ \xi, \zeta \in \mathrm{T}_{X}\Stnp $, the canonical metric reads
\begin{equation}\label{eq:formula_canonical_metric}
    g_{c}(\xi,\zeta) = \trace\!\big( \xi\tr ( I - \tfrac{1}{2}XX\tr) \, \zeta \big).
\end{equation}
The canonical metric induces the \emph{canonical norm}, defined as
\[
    \| \xi \|_{\mathrm{c}} = \sqrt{g_{c}(\xi,\xi)}.
\]
The reader can verify that
\[
   \| \xi \|_{\mathrm{c}}^{2} = \tfrac{1}{2} \| \Omega \|_{\F}^{2} + \| K \|_{\F}^{2}.
\]
The \emph{embedded metric} is the metric inherited by the Stiefel manifold as an embedded submanifold of $ \R^{n} $, i.e., $ g_{e}(\xi,\xi) = \trace( \xi\tr \xi ) $. With the embedded metric, the induced norm is simply the Frobenius norm
\[
    \| \xi \|_{\mathrm{e}} = \sqrt{g_{e}(\xi,\xi)} = \sqrt{\trace( \xi\tr \xi )} \eqqcolon \| \xi \|_{\F}, \qquad \text{and} \qquad   \| \xi \|_{\F}^{2} = \| \Omega \|_{\F}^{2} + \| K \|_{\F}^{2}.
\]
The only difference with respect to the squared canonical norm of $ \xi $ is that in the squared embedded norm of $\xi$ the term $\| \Omega \|_{\F}^{2}$ is not halved. This calculation highlights the fact that in contrast to the embedded norm, the canonical norm only takes into account once the $ \tfrac{1}{2} p (p-1) $ coefficients of $ \Omega $. Indeed, in the remaining part of this paper, we will only use the canonical metric. 

By endowing the Stiefel manifold with the canonical metric, one can derive the following \emph{second-order ordinary differential equation for the geodesic} \protect{\cite[(2.41)]{Edelman:1998}}
\begin{equation}\label{eq:ode_geodesics}
   \ddt{Y} + \dt{Y}\dt{Y}\tr Y + Y \big( (Y\tr\dt{Y})^{2} + \dt{Y}\tr\dt{Y}\big) = 0,
\end{equation}
where $Y \equiv Y(t)$.
An explicit formula for a geodesic that realizes a tangent vector $ \xi $ with base point $ Y_{0} $ is \protect{\cite[(2.42)]{Edelman:1998}}
\begin{equation}\label{eq:closed-form-sol-geodesic}
    Y(t) = Q \expm\!\left( \begin{bmatrix}
        \Omega   &   -K\tr \\
        K        &    O_{n-p}
    \end{bmatrix} t \right)
    \cdot I_{n,p},
\end{equation}
with $ Q = \big[ Y_{0} \ Y_{0\perp} \big] $, $Y_{0\perp}$ being any matrix whose columns span $\mathcal{Y}_{0}^{\perp}=(\mathrm{span}(Y_{0}))^{\perp}$.
If $ t = 1 $, this is precisely the Riemannian exponential on the Stiefel manifold. In this paper, we denote by $ A $ the matrix in the argument of the matrix exponential $ \expm $.

\begin{remark}\label{rmk:Y_0perp}
The matrix $Y_{0\perp}$ does not need to be orthonormal. Indeed, its only requirement is that it has to span $\mathcal{Y}_{0}^{\perp}=(\mathrm{span}(Y_{0}))^{\perp}$, i.e., the orthogonal subspace to $\mathcal{Y}_{0}=\mathrm{span}(Y_{0})$; see~\protect{\cite[Appendix~A.1]{Sutti:2020b}}. For the convenience of our analysis and implementation, we always assume that $Y_{0\perp}$ is orthonormal so that $ Q = \big[ Y_{0} \ Y_{0\perp} \big] $ is an orthogonal matrix.
\end{remark}

\begin{remark}\label{rmk:smaller_formulation}
It can be shown that the endpoint geodesic problem on $\Stnp$ is equivalent to an endpoint geodesic problem on $\mathrm{St}(2p,p)$; see~\protect{\cite{Edelman:1998,Rentmeesters:2013}} and~\protect{\cite[Section~2.3.3, and Appendix A.2]{Sutti:2020b}}.
In the formulation~\eqref{eq:closed-form-sol-geodesic} above, the complexity of computing the matrix exponential is $O(n^3)$, but if $p \ll n $, then the smaller formulation can be used, and its computational cost is only $O(p^3)$. In practice, it makes sense to consider the formulation on $\mathrm{St}(2p,p)$ only if $ p < \tfrac{n}{2} $, which is what we exploit in practice in our algorithm.
\end{remark}

\section{Problem statement}\label{sec:problem_statement}

In this section, we state the problem more formally.
Given two points $Y_0$, $Y_1$ on $\Stnp$ that are sufficiently close to each other, finding the distance between them is equivalent to finding the tangent vector $\xi^{\ast} \in \mathrm{T}_{Y_0}\Stnp$ with the shortest possible length such that \protect{\cite{Lee:2018,boumal2023intromanifolds}}
\[
\Exp_{Y_0}(\xi^{\ast}) = Y_1,
\]
where $\Exp_{Y_0}$ denotes the Riemannian exponential mapping at $Y_0$. 
The solution to this problem is equivalent to the Riemannian logarithm of $Y_{1}$ with base point $Y_{0}$
\[
 \xi^{\ast} = \Log_{Y_0}(Y_1).
\]
We refer the reader to Figure~\ref{fig:single_shooting} for an illustration of the problem statement.

In terms of the differential equation \eqref{eq:ode_geodesics} governing the geodesic, the problem statement may be written as follows:

Find $ \xi^{\ast} \equiv \dt{Y}(0) \in \mathrm{T}_{Y_0}\Stnp$ such that the second-order ODE
\begin{equation}\label{eq:bvp}
\ddt{Y} = - \dt{Y}\dt{Y}\tr Y - Y \big( (Y\tr\dt{Y})^{2} + \dt{Y}\tr\dt{Y}\big), \quad
\text{with boundary conditions} \
\begin{cases}
Y(0)   = Y_{0}, \\
Y(1) = Y_{1},
\end{cases}
\end{equation}
is satisfied. This problem is known as a \emph{boundary value problem} (BVP).

\section{A single shooting method with an approximate Fr\'{e}chet derivative} \label{sec:SSAF}

The single shooting method is a classical numerical scheme for solving boundary value problems. The main idea is to reformulate the BVP as an initial value problem (IVP), guess the initial value of the acceleration, and then solve a nonlinear equation. It basically turns a BVP into a root-finding problem. The zeros of the nonlinear equation can be computed with any root-finding algorithm, but the classical single shooting method typically uses Newton's method.

In this section, we present the details on how to apply the single shooting method to the endpoint geodesic problem on the Stiefel manifold.
We start by recasting the BVP \eqref{eq:bvp} into an IVP.
Let $Z_{1}(t) = Y(t) $, $Z_{2}(t) = \dt{Y}(t)$ denote the geodesic and its derivative, respectively, and let
\[
Z(t) = \begin{pmatrix}
Z_1(t) \\
Z_2(t)
\end{pmatrix}.
\]
We get the initial value problem (we omit the dependence on $t$)
\begin{equation}\label{eq:ivp}
\begin{split}
\dt{Z} = \begin{pmatrix}
\dt{Z}_1 \\
\dt{Z}_2
\end{pmatrix} = \begin{pmatrix}
Z_2 \\
- Z_2 Z_2\tr Z_1 - Z_1 \big( (Z_1\tr Z_2)^{2} + Z_2\tr Z_2 \big)
\end{pmatrix}, \\
\text{with initial conditions} \
Z(0) = \begin{pmatrix}
Z_1(0) \\
Z_2(0)
\end{pmatrix} =
\begin{pmatrix}
Y_0 \\
\xi
\end{pmatrix}.
\end{split}
\end{equation}
Here, $ \xi $ is the unknown such that $ Z_{1}(1) = Y_{1} $.

Solving~\eqref{eq:ivp} typically requires a numerical integration scheme, but here, since we already have the explicit formula \eqref{eq:closed-form-sol-geodesic} for the geodesic $Z_1(t) $, we do not need to integrate the initial value problem \eqref{eq:ivp}.

The explicit formula for $Z_2$ is just the derivative of $Z_1$ with respect to $t$, namely,
\[
Z_2(t) = Q \, \expm\!\left(
\begin{bmatrix}
\Omega   &   -K\tr \\
K        &    O_{n-p}
\end{bmatrix} t \right)
\begin{bmatrix}
\Omega  \\
K     
\end{bmatrix},
\]
where $ \Omega = Y_{0}\tr\xi $ and $ K = Y_{0\perp}\tr \xi $.
Now, let us define the function
\begin{equation}\label{eq:nonlinear-eq}
   F(\xi) = Z_{1}(1,\xi) - Y_1,
\end{equation}
where we emphasize the dependence on $ \xi $.
Roughly speaking, this represents the mismatch between $ Z_{1}(1,\xi) $, i.e., the geodesic at $ t = 1$, and the boundary condition $Y_1$ we wish to enforce.
Our goal is to find $\xi^{\ast}$ such that
\[
   F(\xi^{\ast}) = 0.
\]
As mentioned above, this is a root-finding problem of a nonlinear (matrix) equation, which can be solved by \emph{Newton's method}. To apply Newton's method, we need the Fr\'{e}chet derivative of $Z_{1}(1,\xi)$ in the direction of an increment $\delta\xi$ of $\xi$.

In an earlier version of this work~\protect{\cite{Sutti:2020b,Sutti:2023}}, we used a vectorization approach that allowed us to calculate an explicit analytic expression for the Jacobian matrix involved in this single shooting method. The drawback of this approach is that it is very computationally inefficient due to the dimensions of the operators involved, which grow exponentially with $n$. Here, instead, we work directly with matrix equations and Fr\'{e}chet derivatives.
We draw some inspiration from \protect{\cite{Nguyen:2022}}.

In the remaining part of this section, we first state the algorithm, and then, in Section~\ref{sec:linearization_42}, explain in more detail the linearization of~\eqref{eq:nonlinear-eq} and the approximation of the Fr\'{e}chet derivative involved. In Section~\ref{sec:SS_initial_guess}, we provide a way to construct an initial iterate for the algorithm.

The pseudocode for the single shooting method on the Stiefel manifold is given in Algorithm~\ref{algo:SSAF}. Note at line 8 the additional projection step onto the tangent space to ensure that the updated tangent vector is indeed an element of $ \mathrm{T}_{Y_{0}}\Stnp  $. We did not take any particular care in introducing a line-search technique, although it might be helpful for globalizing the method. Numerical experiments in Section~\ref{sec:numerical_experiments} demonstrate that Algorithm~\ref{algo:SSAF} works very well in practice. For the sake of brevity, we name our algorithm SSAF (an acronym for ``single shooting approximate Fr\'{e}chet''). As a stopping criterion, we typically consider a tolerance on the norm of the residual $ \delta\xi^{(k)} $.

%===========================================================================
% ALGORITHM 1: SINGLE SHOOTING METHOD
%===========================================================================
\begin{algorithm}
   \SetAlgoLined
   Given $Y_0$, $Y_1$\;
   \KwResult{$\xi^{\ast}$ such that $\Exp_{Y_0}(\xi^{\ast}) = Y_1$.}
   Compute the initial guess $ \xi^{(0)} $ (using Algorithm~\ref{algo:initial-guess})\;
   Set $k=0$\;
   \While{a stopping criterion is met}{
      Compute $F^{(k)} = Z_{1}(1,\xi^{(k)}) - Y_1 $\;
      Solve $F^{(k)} + \D Z_{1}\, [\delta\xi^{(k)}] = 0$ for $\delta\xi^{(k)}$\;
      Update $\xi^{(k+1)} \leftarrow \xi^{(k)} + \delta\xi^{(k)}$\;
      Project $ \xi^{(k+1)} $ onto $ \mathrm{T}_{Y_{0}}\Stnp $ using~\eqref{eq:proj_tg_space}: $\xi^{(k+1)} \leftarrow \P_{Y_{0}}\!\big(\xi^{(k+1)}\big)$\;
      $ k = k+1 $\;
   }
   \caption{A single shooting method on the Stiefel manifold with an approximation of the Fr\'{e}chet derivative (SSAF method).}\label{algo:SSAF}
\end{algorithm}
%===========================================================================

\subsection{Linearization of the nonlinear matrix equation~(\ref{eq:nonlinear-eq})} \label{sec:linearization_42}

We recall from~\eqref{eq:tg_space_stiefel} the structure of a tangent vector $\xi \in \mathrm{T}_{Y_{0}}\Stnp $, i.e.,
\[
   \xi = Y_{0} \Omega + Y_{0\perp} K.
\]
From now on, let us denote
\begin{equation}\label{eq:matrix_A}
   A(\xi) = \begin{bmatrix}
             \Omega   &   -K\tr \\
             K        &   O_{n-p}
          \end{bmatrix}
\end{equation}
the matrix in the argument of the exponential appearing in the geodesic equation~\eqref{eq:closed-form-sol-geodesic}. Clearly, $A$ is a function of $ \xi $ because the matrices $ \Omega $ and $ K $ are formed from the tangent vector $ \xi $. Then~\eqref{eq:closed-form-sol-geodesic} at $ t=1 $ can be rewritten as
\[
   Z_1(1,\xi) = Q \expm ( A(\xi) ) \cdot I_{n,p}.
\]
Recall our nonlinear matrix equation~\eqref{eq:nonlinear-eq} that we want to solve for $\xi$.
Newton's method consists in solving successive linearizations~\eqref{eq:nonlinear-eq}, i.e.,
\begin{equation}\label{eq:pert-eq}
   F(\xi + \delta \xi) = Z_1(\xi + \delta \xi) - Y_1 = 0.
\end{equation}

Here, the term $Z_1(\xi + \delta \xi)$ is the expression for the geodesic when applying a small perturbation $ \delta \xi $ to the vector $ \xi $. Applying matrix perturbation theory, we obtain
\[
   Z_1( \xi + \delta \xi ) = Z_1(\xi) +  \D Z_{1}\, [\delta\xi^{(k)}] + o(\Vert \delta \xi \Vert),
\]
namely,
\begin{equation}\label{eq:perturbation-geodesic}
   Z_1( \xi + \delta \xi ) = Z_1(\xi) + Q \, \D \expm(A(\xi))\big[\D \! A(\xi)[\delta \xi]\big] \cdot I_{n,p} + o(\Vert \delta \xi \Vert),
\end{equation}
where $ \D \expm(A(\xi))\big[\D \! A(\xi)[\delta \xi]\big] $ denotes the Fr\'{e}chet derivative of the matrix exponential at $A(\xi)$ in the direction of $ \D \! A(\xi)[\delta \xi] $. $ \D \! A(\xi)[\delta \xi] $ itself denotes the Fr\'{e}chet derivative of $A(\xi)$ in the direction of $ \delta \xi $.

In contrast to what was proposed in~\protect{\cite{Sutti:2020b,Sutti:2023}}, here we do not vectorize the equation, and we work directly with matrices and Fr\'{e}chet derivatives; moreover, we do not compute the exact Fr\'{e}chet derivative, but we approximate it by a truncated expansion.

Inserting~\eqref{eq:perturbation-geodesic} into~\eqref{eq:pert-eq} and neglecting the higher-order terms in $\delta\xi$, we obtain the matrix equation
\begin{equation}\label{eq:4pt8}
   Z_1(\xi) + Q \, \D \expm(A(\xi))\big[\D \! A(\xi)[\delta \xi]\big] \cdot
            I_{n,p}
            - Y_{1} = 0.
\end{equation}
We now need to tackle the term $ \D \expm(A(\xi))\big[\D \! A(\xi)[\delta \xi]\big] $, which involves a chain rule with two Fr\'{e}chet derivatives: the term $\D \! A(\xi)[\delta \xi]$ is the Fr\'{e}chet derivative of $ A(\xi) $ in the direction of $ \delta\xi $; and $\D \expm(A(\xi))\big[\D \! A(\xi)[\delta \xi]$, that is the Fr\'{e}chet derivative of the matrix exponential $\expm(A(\xi))$ in the direction of $\D \! A(\xi)[\delta \xi]$.

First, the perturbation of $A(\xi)$ with a $\delta\xi$ gives
\[
   A(\xi+\delta \xi ) = A(\xi) +  \D \! A(\xi)[\delta \xi].
\]
where
\begin{equation}
   A(\xi) =
   \begin{bmatrix}
      \Omega   &   -K\tr   \\
      K        &   O_{(n-p)}
   \end{bmatrix}, \quad \text{and} \quad
   \D \! A(\xi)[ \delta \xi ] =
   \begin{bmatrix}
      \delta \Omega   &   -\delta K\tr   \\
      \delta K        &   O_{(n-p)}
   \end{bmatrix}.
\end{equation}
Since $A$ is linear in $\xi$, the above expansion is exact.

Secondly, the perturbation of the matrix exponential by a matrix $ E \in \R^{n \times n} $ is
\[
   \expm( A + E ) = \expm(A) + \D \expm(A)[E] + o(\Vert E \Vert),
\]
where $ \D \expm(A) [E] $ is the Fr\'{e}chet derivative of the matrix exponential at $A$ in the direction of $ E $. In general, there are many ways to compute $\expm$ and $\D \expm(A)[E]$, and thus also many ways to approximate these quantities; see~\protect{\cite[Chapter~10]{Higham:2008}}. Here, we consider the Taylor series of $ e^{A+E} - e^{A} $, from which we obtain the following representation for the Fr\'{e}chet derivative of the matrix exponential~\protect{\cite[Section~10.2]{Higham:2008}}
\[
   \D \expm(A)[E] = E + \frac{AE + EA}{2} + \frac{A^{2}E + AEA + EA^{2}}{3!} + \cdots .
\]
We then consider an approximation to $\D \expm(A)[E]$ by keeping only the first two terms in the expansion, i.e.,
\[
   \D \expm(A)[E] \approx E + \frac{AE + EA}{2}.
\]
This formula can be used to approximate $\D \expm(A(\xi))\big[\D \! A(\xi)[\delta \xi]\big]$ in~\eqref{eq:4pt8}, resulting in
\[
   Q \cdot \left( \D \! A(\xi)[\delta \xi] + \tfrac{1}{2} \! \left( A \cdot \D \! A(\xi)[\delta \xi] + \D \! A(\xi)[\delta \xi] \cdot A \right) \right) \cdot I_{n,p}
            = Y_{1} - Z_{1},
\]
This is now a linear matrix equation to be solved for $\delta\xi$. In practice, we work with the factors $ \delta\Omega $ and $ \delta K $ of $ \delta\xi $. After a few algebraic manipulations, detailed in Appendix~\ref{app:approx_frechet}, we obtain a (small-sized) Sylvester equation which can be efficiently solved with MATLAB's command \texttt{lyap} to obtain the update $\delta \Omega$. Then the update $\delta K$ can be found from $\delta \Omega$. Let the current iteration be indexed by $k$; then the tangent vector is updated as
\[
   \xi^{(k+1)} = \xi^{(k)} + Q \cdot
   \begin{bmatrix}
      \delta \Omega^{(k)} \\
      \delta K^{(k)}   
   \end{bmatrix}.
\]

\subsection{The initial guess}\label{sec:SS_initial_guess}

This section outlines our approach to initializing our single shooting method, which involves choosing an initial guess $\xi^{(0)}$ that is close enough to $\xi^{\ast}$.
It is well known that Newton's method exhibits only local convergence properties, which means that the method requires a sufficiently good initial guess to converge. Shortcomings of Newton's method are very well described in~\protect{\cite[Section~11.1]{Nocedal:2006}}. It is possible to modify Newton's method and enhance it in various ways to get around most of these problems.

As the single shooting method underlies Newton's method, selecting a ``good enough'' initial guess is crucial. Although in this work we actually consider an approximation of the Fr\'{e}chet derivative and not the exact one, we can still use the following construction of the initial guess that was used in~\protect{\cite{Sutti:2020b,Sutti:2023}}. This construction is closely related to the ``first shot'' in Bryner's method~\protect{\cite[Alg.~1]{Bryner:2017}}.

Concretely, we use a first-order approximation of the matrix exponential $ \expm(A) \approx I + A $ appearing in~\eqref{eq:nonlinear-eq} and solve for $\xi$. 
This yields a first-order approximation $\bar{\xi}$ to the solution $\xi^{\ast}$ as
\[
   \bar{\xi} = Y_1 - Y_0.
\]
Since, in general, this is no longer an element of the tangent space, we need to project it onto $ \mathrm{T}_{Y_0}\Stnp $ to obtain a tangent vector. We expect the tangent vector so obtained to be a satisfactory initial approximation to the sought tangent vector $\xi^{\ast}$.

\begin{figure}[htbp]
   \centering
   \includegraphics[width=0.55\columnwidth]{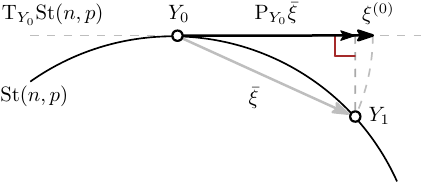}
   \caption{Initial guess for the single shooting method on the Stiefel manifold.}
   \label{fig:SS_initial_guess}
\end{figure}

Using~\eqref{eq:proj_tg_space}, the projection of $\bar{\xi}$ onto the tangent space at $Y_0$ is
\[
   \P_{Y_{0}} \! \bar{\xi} = Y_0 \,\mathrm{skew}\big(Y_0\tr(Y_1 - Y_0)\big) + (I_n - Y_0 Y_0\tr)(Y_1 - Y_0) = Y_1 - Y_0 \, \mathrm{sym}(Y_0\tr Y_1).
\]
To get $\xi^{(0)}$, we rescale this vector so that its norm is equal to the norm of $\bar{\xi}$, i.e.,
\[
   \xi^{(0)} = \frac{\left\Vert \bar{\xi} \right\Vert}{\left\Vert \P_{Y_{0}} \! \bar{\xi} \right\Vert} \, \P_{Y_{0}} \! \bar{\xi}.
\]
This procedure is summarized in Algorithm~\ref{algo:initial-guess} and illustrated in Figure~\ref{fig:SS_initial_guess}.

%=========================================
% FIGURE: INITIAL GUESS SINGLE SHOOTING
%=========================================
%=========================================

%===========================================================================
% ALGORITHM 2: 
%===========================================================================
\begin{algorithm}
   \SetAlgoLined
   Given $Y_0$, $Y_1$\;
   Compute $ \bar{\xi} = Y_1 - Y_0 $\;
   Compute $ \P_{Y_{0}} \! \bar{\xi} = Y_1 - Y_0 \, \mathrm{sym}(Y_0\tr Y_1) $\;
   Compute $ \xi^{(0)} = \frac{\left\Vert \bar{\xi} \right\Vert}{\left\Vert \P_{Y_{0}} \! \bar{\xi} \right\Vert} \, \P_{Y_{0}} \! \bar{\xi} $\;
   Return  $ \xi^{(0)} $.
   \caption{Initial guess for the single shooting method on the Stiefel manifold.}\label{algo:initial-guess}
\end{algorithm}
%===========================================================================

\section{Numerical experiments and comparisons with other methods}\label{sec:numerical_experiments}

In this section, we present some numerical experiments about the single shooting method, and we report on the convergence behavior. The code was implemented in MATLAB and is freely available on the repository \href{https://github.com/MarcoSutti/SSAF_2024_repo}{https://github.com/MarcoSutti/SSAF\_2024\_repo}. The method of Bryner~\protect{\cite[Alg.~1]{Bryner:2017}} was implemented according to the pseudocode provided in~\protect{\cite{Bryner:2017}}, while the method of Zimmermann~\protect{\cite[Alg.~1]{Bryner:2017}} was directly taken from Appendix C of~\protect{\cite{Zimmermann:2017}}.
We conducted our experiments on a laptop Lenovo ThinkPad T460s with Ubuntu 23.10 LTS and MATLAB R2022a installed, with Intel Core i7-6600 CPU, 20GB RAM, and Mesa Intel HD Graphics 520. On this machine, the matrix exponential \texttt{expm} of a unit norm skew-symmetric matrix in $\R^{1000 \times 1000}$ is computed in 0.45 seconds (time averaged over 100 runs)\footnote{We give this reference for the computational cost of the matrix exponential \texttt{expm} following~\protect{\cite[Section~5.2]{Nguyen:2022}}.}. In all the tables, our algorithm is named SSAF.

Table~\ref{tab:n1} compares our new SSAF method with the earlier version using the exact Jacobian matrix~\protect{\cite[Alg.~1]{Sutti:2020b}}. The efficiency of our new SSAF method compared to the old version is striking. Typically, cases with a small $p$ require more iterations for SSAF, while for~\protect{\cite[Alg.~1]{Sutti:2020b}} the number of iterations remains constant w.r.t. $p$. Still, the superior efficiency of SSAF by far compensates for this (desirable) feature.
From Table~\ref{tab:n1}, it is also evident that the single shooting method with the exact Jacobian matrix~\protect{\cite[Alg.~1]{Sutti:2020b}} scales very badly with $p$, and it becomes prohibitively expensive as $p$ grows. The long dashes ``---'' in the table indicate that the single shooting method with the exact Jacobian matrix stopped due to memory overflow.

%=================================================================================
% TABLE 1: COMPARISON SS WITH WITH EXACT JACOBIAN VS SS WITH APPROXIMATE FRECHET
%=================================================================================
\begin{table}[htbp]
   \caption{Comparisons for the single shooting method with approximate Fr\'{e}chet derivative versus the single shooting method with the exact Jacobian~\protect{\cite[Alg.~1]{Sutti:2020b}}, on $\mathrm{St}(1000,p)$, for doubling values of $p$, for a prescribed $ d(X,Y) = 0.5\,\pi $. Results are averaged over 10 random runs. Stopping tolerance is $10^{-5}$.}
   \label{tab:n1}
   \begin{center}
      \begin{tabular}{@{} l *{4}{c}}
         \multirow{2}{*}{$p$}   &        \mcii{Avg. comput. time (s)}         &      \mcii{Avg. no. of iterations}      \\
                                       \cmidrule(lr){2-3}              \cmidrule(lr){4-5}
                                &   \protect{\cite[Alg.~1]{Sutti:2020b}}  &   SSAF     &   \protect{\cite[Alg.~1]{Sutti:2020b}}  &   SSAF   \\
         \midrule
                        10      &       0.08   &     0.00256     &         5  &     6.80     \\
                        20      &       3.07   &     0.00390     &         5  &     5.00     \\
                        40      &     182.69   &     0.00998     &         5  &     5.00     \\
                        80      &       ---    &     0.03388     &     ---    &     4.00     \\
                       160      &       ---    &     0.12734     &     ---    &     4.00     \\
                       320      &       ---    &     0.90012     &     ---    &     4.00     \\                    
                       640      &       ---    &     2.98705     &     ---    &     4.00
      \end{tabular}
   \end{center}
\end{table}
%=================================================================================

\subsection{Comparisons with other state-of-the-art methods}\label{sec:comparisons}

In this section, we demonstrate that the proposed algorithm is competitive with other state-of-the-art methods.
The body of literature and available methods have been increasing recently, especially during the last five years, and it would be hard to compare our proposed algorithm to all the existing methods. Hence, the comparisons in this section are not meant to be exhaustive. We aim to show that our proposed algorithm is competitive with respect to only a limited subset of state-of-the-art algorithms that can be found in the literature.
Specifically, in this section we compare our method to the ``shooting'' method of Bryner~\protect{\cite[Alg.~1]{Bryner:2017}}, the matrix algebraic approach of Zimmermann~\protect{\cite[Alg.~1]{Zimmermann:2017}}, and the optimization methods of Nguyen~\protect{\cite{Nguyen:2022}}.

Table~\ref{tab:n2} uses the same test cases as those considered in~\protect{\cite[Table~2]{Nguyen:2022}}, namely $\mathrm{St}(1500,p)$, for large values of $p$ and for a prescribed $ d(X,Y) = 0.5\,\pi $\footnote{This is the prescribed distance used for the numerical experiments reported in~\protect{\cite[Table~2]{Nguyen:2022}}, although not explicitly written in that paper. From a private communication with Nguyen, February 2024.}. In all the test cases considered, our SSAF algorithm is superior in terms of computation time compared to the other methods.

%=====================================================
% TABLE 2: SIMILAR TO EXPERIMENTS OF DU NGUYEN, 2022
%=====================================================
\begin{table}[htbp]
   \caption{Comparisons on $\mathrm{St}(1500,p)$ with large values of $p$, for a prescribed $ d(X,Y) = 0.5\,\pi $. Results are averaged over 10 pairs of randomly generated endpoints on $\mathrm{St}(1500,p)$. Stopping tolerance is $10^{-5}$.}
   \label{tab:n2}
   \begin{center}
      \begin{tabular}{@{} l *{6}{c}}
         \multirow{2}{*}{$p$}   &        \mciii{Avg. comput. time (s)}         &      \mciii{Avg. no. of iterations}      \\
                                       \cmidrule(lr){2-4}              \cmidrule(lr){5-7}
                                &   \protect{\cite[Alg.~1]{Bryner:2017}}   &   \protect{\cite[Alg.~1]{Zimmermann:2017}}   &    SSAF     &   \protect{\cite[Alg.~1]{Bryner:2017}}   &   \protect{\cite[Alg.~1]{Zimmermann:2017}}   &      SSAF     \\
         \midrule
                        500     &   12.09   &      3.30     &  1.89   &     3.00     &       2.00       &     4.00    \\
                        700     &   31.27   &      8.21     &  4.56   &     3.00     &       2.00       &     4.00    \\
                       1000     &   77.37   &     20.39     &  8.82   &     3.00     &       2.00       &     4.00
      \end{tabular}
   \end{center}
\end{table}
%=====================================================

Here, we do not perform a direct comparison with the methods of~\protect{\cite{Nguyen:2022}}, in particular, the results in~\protect{\cite[Table~2]{Nguyen:2022}. However, we can point out that the results in~\protect{\cite[Table~2]{Nguyen:2022}} are obtained on a machine that computes the matrix exponential \texttt{expm} of a unit norm skew-symmetric matrix in $\R^{1000 \times 1000}$ is computed in 0.6 seconds. In contrast, on our machine, as mentioned at the beginning of this section, this same reference quantity is 0.45 seconds. This suggests that if it were to be run on the same machine, the methods of~\protect{\cite{Nguyen:2022}} would be slightly faster but still slower than our SSAF method and also than~\protect{\cite[Alg.~1]{Zimmermann:2017}}.

However, when $p$ is doubled from 500 to 1000, the methods of~\protect{\cite[Table~2]{Nguyen:2022}} seem to scale better than all the other algorithms considered here. 
We emphasize that the test case corresponding to the last row in Table~\ref{tab:n2} has dimensions $n=1500$ and $p=1000$, i.e., this is not a case in which $ p \leq n/2 $. Therefore, the standard implementations of \protect{\cite[Alg.~1]{Bryner:2017}} and~\protect{\cite[Alg.~1]{Zimmermann:2017}} are not designed to be efficient in this case.
From the numerical results reported in~\protect{\cite[Table~2]{Nguyen:2022}}, it seems that when $p$ is doubled from 500 to 1000, computation time is multiplied by a factor of approximately 3.4 for the gradient descent method, and 3.6 for the L-BFGS method. In contrast, the same factors computed from Table~\ref{tab:n2} for~\protect{\cite[Alg.~1]{Bryner:2017}}, \protect{\cite[Alg.~1]{Zimmermann:2017}}, and our SSAF method are 6.4, 6.2, and 4.7, respectively. This suggests that Nguyen's methods might be more effective for problems with larger values of $p$. Yet our SSAF remains the most competitive among the other algorithms considered here because it shows to have a factor of 4.7, in contrast to 6.4 for~\protect{\cite[Alg.~1]{Bryner:2017}} and 6.2 for~\protect{\cite[Alg.~1]{Zimmermann:2017}}.

Table~\ref{tab:n3} considers the same three test cases as in~\protect{\cite[Table~1]{Zimmermann:2022}}, with stopping tolerance $ \tau = 10^{-10} $. Although the average number of iterations for our SSAF method is much higher than that of the other algorithms, our method remains competitive in terms of computation time.

%===========================================================
% TABLE 3: TEST CASES IN TABLE 1 OF ZIMMERMANN HUPER 2022
%===========================================================
\begin{table}[htbp]
   \caption{Comparisons for the three test cases in~\protect{\cite[Table~1]{Zimmermann:2022}}, with stopping tolerance $ \tau = 10^{-10} $. Results are averaged over 10 experiments.}
   \label{tab:n3}
   \begin{center}
      \begin{tabular}{@{} *{6}{c}}
         \mciii{Avg. comput. time (s)}         &      \mciii{Avg. no. of iterations}      \\
         \cmidrule(lr){1-3}              \cmidrule(lr){4-6}
         \protect{\cite[Alg.~1]{Bryner:2017}}  &   \protect{\cite[Alg.~1]{Zimmermann:2017}}   &    SSAF     &   \protect{\cite[Alg.~1]{Bryner:2017}}   &   \protect{\cite[Alg.~1]{Zimmermann:2017}}   &         SSAF     \\
         \midrule
         \multicolumn{6}{c}{Test Case 1 : $\mathrm{St}(2000,500)$, for a prescribed $ d(X,Y) = 5\pi $.} \\
         \midrule
            109.65  &       15.35      &    19.66   &     14.0     &        12.0      &        28.2      \\
         \midrule
         \multicolumn{6}{c}{Test Case 2: $\mathrm{St}(120,30)$, for a prescribed $ d(X,Y) = \pi $.} \\
         \midrule
            0.11242  &     0.04157     &  0.02556   &     10.90    &        9.30      &        19.60     \\
         \midrule
         \multicolumn{6}{c}{Test Case 3: $\mathrm{St}(12,3)$, for a prescribed $ d(X,Y) = 0.95\pi $.} \\
         \midrule
            0.08885  &     0.06563     &  0.03948   &     40.80    &       86.30     &       167.00
      \end{tabular}
   \end{center}
\end{table}
%===========================================================

Tables~\ref{tab:n4} and~\ref{tab:n5} below try to reproduce the data from the left and middle panels, respectively, of~\protect{\cite[Figure~4]{Bryner:2017}}, while at the same time comparing with the method of~\protect{\cite{Zimmermann:2017}} and with our SSAF method. We adopt the same parameters as in~\protect{\cite{Bryner:2017}}, namely, $T=20$, and stopping tolerance $\tau = 10^{-3}$.
We say that we ``try to reproduce'' since it seems that Bryner did not prescribe a distance, or at least this is not explicitly stated; hence, we fix it here to $ d(X,Y) = 0.5\,\pi $.

%=====================================================
% TABLE 4: COMPARISON WITH OTHER METHODS, DOUBLING n.
%          Reproduction of left panel of Figure 4 in Bryner, 2017.
%=====================================================
\begin{table}[htbp]
   \caption{Comparisons on $\mathrm{St}(n,2)$, for doubling values of $ n $, for a prescribed $ d(X,Y) = 0.5\,\pi $. $T=20$, tolerance $ \tau = 10^{-3} $. Results are averaged over 100 experiments.}
   \label{tab:n4}
   \begin{center}
      \begin{tabular}{@{} l *{6}{c}}
         \multirow{2}{*}{$n$}   &        \mciii{Avg. comput. time (s)}       &      \mciii{Avg. no. of iterations}      \\
                                       \cmidrule(lr){2-4}              \cmidrule(lr){5-7}
                                &  \protect{\cite[Alg.~1]{Bryner:2017}}  &   \protect{\cite[Alg.~1]{Zimmermann:2017}}   &    SSAF      &   \protect{\cite[Alg.~1]{Bryner:2017}}  &   \protect{\cite[Alg.~1]{Zimmermann:2017}}   &    SSAF     \\
         \midrule
                        10      &   0.00400  &      0.00091     &   0.00080   &    4.08     &       3.73       &    7.77    \\
                        20      &   0.00367  &      0.00093     &   0.00091   &    3.85     &       3.87       &    7.35    \\
                        40      &   0.00337  &      0.00095     &   0.00075   &    3.49     &       3.61       &    6.96    \\
                        80      &   0.00312  &      0.00101     &   0.00081   &    3.30     &       3.61       &    6.90    \\
                       160      &   0.00310  &      0.00105     &   0.00086   &    3.15     &       3.42       &    6.86    \\
                       320      &   0.00328  &      0.00107     &   0.00096   &    3.02     &       3.08       &    6.86    \\
                       640      &   0.00371  &      0.00105     &   0.00091   &    3.00     &       3.02       &    6.89    \\
                    1\,280      &   0.00543  &      0.00104     &   0.00100   &    3.00     &       2.72       &    6.87    \\
                    2\,560      &   0.00856  &      0.00135     &   0.00121   &    3.00     &       2.47       &    6.87    \\
                    5\,120      &   0.01056  &      0.00131     &   0.00132   &    3.00     &       2.34       &    6.93    \\
                   10\,240      &   0.01596  &      0.00144     &   0.00141   &    3.00     &       2.12       &    6.97
      \end{tabular}
   \end{center}
\end{table}
%=====================================================

From Table~\ref{tab:n5}, it appears that the endpoint geodesic problem seems to get easier for large $p$, since for all the three methods considered, the average number of iterations decreases for increasing $n$. In other words, as the ratio $ p/n \to 1 $, solving the endpoint geodesic problem requires fewer iterations. This is consistent with a similar observation in~\protect{\cite[Section~5.2]{Nguyen:2022}}, and with~\protect{\cite[Table~5.1]{Zimmermann:2017}}.}

%=====================================================
% TABLE 5: COMPARISON WITH OTHER METHODS, DOUBLING n.
%          Reproduction of the middle panel of Figure 4 in Bryner, 2017.
%=====================================================
\begin{table}[htbp]
   \caption{Comparisons on $\mathrm{St}(500,p)$, for doubling values of $ p $, for a prescribed $ d(X,Y) = 0.5\,\pi $. $T=20$, tolerance $ \tau = 10^{-3} $. Results are averaged over 100 experiments.}
   \label{tab:n5}
   \begin{center}
      \begin{tabular}{@{} l *{6}{c}}
         \multirow{2}{*}{$p$}   &        \mciii{Avg. comput. time (s)}       &      \mciii{Avg. no. of iterations}      \\
                                       \cmidrule(lr){2-4}              \cmidrule(lr){5-7}
                                &  \protect{\cite[Alg.~1]{Bryner:2017}}  &   \protect{\cite[Alg.~1]{Zimmermann:2017}}   &    SSAF      &   \protect{\cite[Alg.~1]{Bryner:2017}}  &   \protect{\cite[Alg.~1]{Zimmermann:2017}}   &    SSAF     \\
         \midrule
                         2      &   0.00353  &   0.00103   &   0.00086   &    3.01     &     2.95    &    6.78   \\
                         4      &   0.00533  &   0.00156   &   0.00128   &    3.00     &     2.81    &    5.28   \\
                         8      &   0.00711  &   0.00182   &   0.00115   &    3.00     &     2.00    &    4.08   \\
                        16      &   0.01173  &   0.00369   &   0.00173   &    3.00     &     2.00    &    4.00   \\
                        32      &   0.02912  &   0.01354   &   0.00453   &    3.00     &     2.00    &    4.00   \\
                        64      &   0.08762  &   0.03582   &   0.01150   &    3.00     &     2.00    &    3.00   \\           
                       128      &   0.40437  &   0.10052   &   0.05657   &    3.00     &     1.00    &    3.00   \\
                       256      &   1.94025  &   0.47720   &   0.25847   &    3.00     &     1.00    &    3.00 
      \end{tabular}
   \end{center}
\end{table}
%=====================================================

\begin{figure}[htbp]
  \centering
  \includegraphics[width=\textwidth]{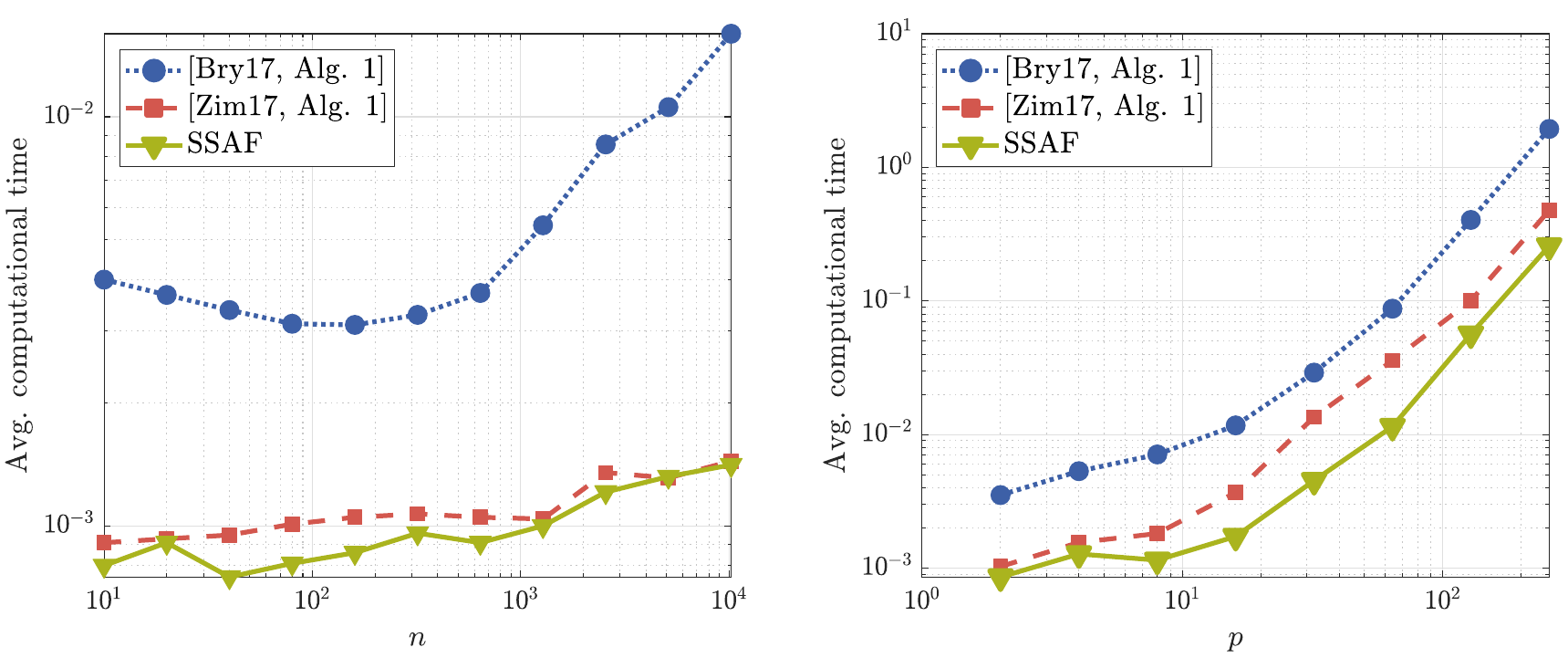}
  \caption{Average computation times for Bryner's shooting method~\protect{\cite[Alg.~1]{Bryner:2017}}, Zimmermann's matrix algebraic algorithm~\protect{\cite[Alg.~1]{Zimmermann:2017}}, and our SSAF method on $\Stnp$. Left panel: plot corresponding to Table~\ref{tab:n5}. Right panel: plot corresponding to Table~\ref{tab:n6}.}\label{fig:Comparison_avg_comput_times}
\end{figure}

Using the data from Tables~\ref{tab:n4} and~\ref{tab:n5}, Figure~\ref{fig:Comparison_avg_comput_times} tries to reproduce the left and middle panels of~\protect{\cite[Figure~4]{Bryner:2017}}, while at the same time comparing with Zimmermann's algorithm~\protect{\cite[Alg.~1]{Zimmermann:2017}} and our SSAF method.

%=========================================
% MATLAB PLOT: FIGURE 4
%=========================================
%=========================================

We emphasize that Bryner did not use the smaller formulation on $\mathrm{St}(2p,p)$ when $p<n/2$ (see Remark~\ref{rmk:smaller_formulation} above), which makes its algorithm's complexity $\cO(Tnp^{2})$; see~\protect{\cite[Section~5.2]{Nguyen:2022}}. The other algorithms considered here (\protect{\cite{Nguyen:2022}}, \protect{\cite[Alg.~1]{Zimmermann:2017}}, and our new SSAF algorithm) all make use of the smaller formulation on $\mathrm{St}(2p,p)$ when possible and hence they are essentially $\cO(p^{3})$, with our SSAF being comparable or even superior than~\protect{\cite[Alg.~1]{Zimmermann:2017}} in terms of average computation time; see Tables~\ref{tab:n4}, \ref{tab:n5}, and Figure~\ref{fig:Comparison_avg_comput_times}.

As a last numerical experiment, we consider a larger value of $n$, namely $n = 1000$, and larger doubling values of $p$. Results are reported in Table~\ref{tab:n6} below, demonstrating again the competitiveness of our SSAF method in terms of both average computation time and number of iterations w.r.t. the existing algorithms considered here.

%=====================================================
% TABLE 6: COMPARISON WITH OTHER METHODS, DOUBLING p
%=====================================================
\begin{table}[htbp]
   \caption{Comparisons on $\mathrm{St}(1000,p)$, for doubling values of $ p $, for a prescribed $ d(X,Y) = 0.5\,\pi $. $T=20$, tolerance $ \tau = 10^{-5} $.  Results are averaged over 100 experiments.}
   \label{tab:n6}
   \begin{center}
      \begin{tabular}{@{} l *{6}{c}}
         \multirow{2}{*}{$p$}   &        \mciii{Avg. comput. time (s)}         &      \mciii{Avg. no. of iterations}      \\
                                       \cmidrule(lr){2-4}              \cmidrule(lr){5-7}
                                &   \protect{\cite[Alg.~1]{Bryner:2017}}  &   \protect{\cite[Alg.~1]{Zimmermann:2017}}   &    SSAF     &   \protect{\cite[Alg.~1]{Bryner:2017}}   &   \protect{\cite[Alg.~1]{Zimmermann:2017}}   &         SSAF     \\
         \midrule
                        20      &    0.03897  &      0.00641     &  0.00391   &     4.00     &        3.00      &        5.02     \\
                        40      &    0.09512  &      0.02957     &  0.01284   &     3.00     &        3.00      &        5.00     \\
                        80      &    0.25528  &      0.08044     &  0.03969   &     3.00     &        2.00      &        4.00     \\
                       160      &    0.76246  &      0.24119     &  0.13763   &     3.00     &        2.00      &        4.00     \\
                       320      &    3.99810  &      1.07286     &  0.64483   &     3.00     &        2.00      &        4.00     \\                       
                       640      &   23.36386  &      5.62897     &  2.80133   &     3.00     &        2.00      &        4.00
      \end{tabular}
   \end{center}
\end{table}
%=====================================================

\section{Conclusions and outlook} \label{sec:conclusions}

In this work, we studied the shooting method, a classical numerical algorithm for solving boundary value problems, to compute the distance between two given points on the Stiefel manifold under the canonical method. We provided a shooting method for calculating geodesics on the Stiefel manifold in the sense of classical shooting methods for solving boundary value problems. The main feature of our algorithm is that we provide an approximate formula for the Fr\'{e}chet derivative of the geodesic involved in our shooting method.
Numerical experiments demonstrated our algorithm's performance and accuracy. We compared our algorithm to some state-of-the-art methods and showed that it is competitive with existing algorithms.

As a future outlook, an analysis of the proposed algorithm would be desirable. Moreover, we may use the knowledge gained in this work to develop a computationally cheaper algorithm. 
Another promising research direction is exploring the connection between shooting algorithms for computing geodesics and domain decomposition methods. Future studies will focus on these topics.

\section*{Acknowledgments}

The author is grateful to Bart Vandereycken for his guidance during the author's Ph.D. thesis.
Part of this work was started during the author's Ph.D. thesis at the University of Geneva, SNSF fund number 163212\footnote{SNSF webpage: \href{https://data.snf.ch/grants/grant/163212}{https://data.snf.ch/grants/grant/163212}}. It was completed during the author's postdoctoral fellowship at the National Center for Theoretical Sciences in Taiwan (R.O.C.) under the NSTC grant 112-2124-M-002-009-.

The author would also like to thank the anonymous referee for the many valuable comments regarding an earlier version of this paper, which considerably help to improve the manuscript.

\appendix

\section{Approximation of the Fr\'{e}chet derivative of the matrix exponential}\label{app:approx_frechet}

We recall that the Fr\'{e}chet derivative of a matrix function $f \colon \C^{n\times n} \to \C^{n\times n}$ at $X \in \C^{n\times n}$ is the unique linear function $ \D\! f(X) [\cdot] $ of the matrix $E \in \C^{n\times n}$, that satisfies
\begin{equation}\label{eq:1}
   f(X+E) - f(X) - \D\! f(X)[E] = o(\| E \|).
\end{equation}
The mapping itself is denoted by either $ \D\! f(X)[\cdot] $ or $ \D\! f(X) $, while the value of the mapping for direction $E$ (i.e., the directional derivative) is denoted by $\D\! f(X)[E]$.

From~\eqref{eq:perturbation-geodesic}, we have the matrix equation
\begin{equation}\label{eq:C1}
   Z_1(\xi) + Q \, \D \expm(A(\xi))\big[\D \! A(\xi)[\delta \xi]\big]
            I_{n,p}
            - Y_{1} = 0,
\end{equation}
The Fr\'{e}chet derivative of the matrix exponential is defined through the integral~\protect{\cite[(10.15)]{Higham:2008}}
\[
   \D \expm(A)[E] \coloneqq \int_{0}^{1} e^{A(1-s)} E \, e^{As} \, \mathrm{d}s.
\]
We also have the following formula, from the Taylor series of $ e^{A+E} - e^{A} $~\protect{\cite[Section~10.2]{Higham:2008}}
\[
   \D \expm(A)[E] = E + \frac{AE + EA}{2} + \frac{A^{2}E + AEA + EA^{2}}{3!} + \cdots .
\]
As mentioned in Section~\ref{sec:linearization_42}, here we consider an approximation to $\D \expm(A)[E]$ by keeping only the first two terms in the expansion, i.e.,
\[
   \D \expm(A)[E] \approx E + \frac{AE + EA}{2}.
\]
This truncated expansion can be used to approximate $\D \expm(A(\xi))\big[\D \! A(\xi)[\delta \xi]\big]$ in~\eqref{eq:C1}, which yields
\[
   Q \cdot \left( \D \! A(\xi)[\delta \xi] + \tfrac{1}{2} \! \left( A \cdot \D \! A(\xi)[\delta \xi] + \D \! A(\xi)[\delta \xi] \cdot A \right) \right) \cdot I_{n,p}
            = Y_{1} - Z_{1},
\]
Left-multiplying the last equation by $ Q\tr $, we get
\[
   \left( \D \! A(\xi)[\delta \xi] + \tfrac{1}{2} \! \left( A \cdot \D \! A(\xi)[\delta \xi] + \D \! A(\xi)[\delta \xi] \cdot A \right) \right) \cdot I_{n,p}
            = Q\tr \! \left(Y_{1} - Z_{1}\right),
\]
We emphasize that $Q\tr \! Z_{1} = \expm(A) \, I_{n,p}$, and the other term $Q\tr\! Y_{1}$ does not depend on $\xi$; hence, in the practical implementation of the algorithm, we compute this quantity only once. Continuing with the manipulations, we obtain
\[
   \begin{bmatrix}
      \delta \Omega \\
      \delta K
   \end{bmatrix}
   + \frac{1}{2}
   \left(
   \begin{bmatrix}
      \Omega   &   -K\tr   \\
      K        &   O
   \end{bmatrix}
   \begin{bmatrix}
      \delta \Omega   \\
      \delta K  
   \end{bmatrix} +
   \begin{bmatrix}
      \delta \Omega   &   -\delta K\tr   \\
      \delta K        &   O
   \end{bmatrix}
   \begin{bmatrix}
      \Omega  \\
      K 
   \end{bmatrix}
   \right) = Q\tr \! \left(Y_{1} - Z_{1}\right),
\]
\[
   \begin{bmatrix}
      \delta \Omega \\
      \delta K
   \end{bmatrix}
   + \frac{1}{2}
   \left(
   \begin{bmatrix}
      \Omega \delta \Omega  -K\tr \! \delta K    \\
      K \delta \Omega
   \end{bmatrix} +
   \begin{bmatrix}
      \delta \Omega \Omega -\delta K\tr \! K \\
      \delta K \Omega 
   \end{bmatrix}
   \right) = Q\tr \! \left(Y_{1} - Z_{1}\right),
\]
\[
   \begin{bmatrix}
      \delta \Omega + \tfrac{1}{2} \Omega \delta \Omega - \tfrac{1}{2} K\tr \! \delta K + \tfrac{1}{2} \delta \Omega \Omega - \tfrac{1}{2} \delta K\tr \! K \\[6pt]
      \delta K + \tfrac{1}{2} K \delta \Omega + \tfrac{1}{2} \delta K \Omega
   \end{bmatrix}
   = Q\tr \! \left(Y_{1} - Z_{1}\right),
\]
\begin{equation}\label{eq:C3}
   \begin{bmatrix}
      \delta \Omega + \tfrac{1}{2} [\Omega, \delta \Omega ] - \tfrac{1}{2} \left( K\tr \! \delta K + \delta K\tr \! K \right) \\[6pt]
      \delta K + \tfrac{1}{2} K \delta \Omega + \tfrac{1}{2} \delta K \Omega
   \end{bmatrix}
   =
   \begin{bmatrix}
      \left. Q\tr \! \left(Y_{1} - Z_{1}\right) \right|_{[1:p, \colon]} \\
      \left. Q\tr \! \left(Y_{1} - Z_{1}\right) \right|_{[p+1:n, \colon]}
   \end{bmatrix}
   \eqqcolon
   \begin{bmatrix}
      W \\
      N
   \end{bmatrix}
   .
\end{equation}
From the second matrix equation, we have
\begin{equation}\label{eq:delta_K}
   \delta K \left( I_{p} + \tfrac{1}{2} \Omega \right) = N - \tfrac{1}{2} K \delta \Omega,
\end{equation}
and we approximate $\left( I_{p} + \tfrac{1}{2} \Omega \right)$ with $ I_{p} $, i.e.,
\[
   \delta K = N - \tfrac{1}{2} K \delta \Omega.
\]
Now we insert this last equation into the first matrix equation in~\eqref{eq:C3} to solve for $ \delta \Omega $, giving
\[
   \delta \Omega + \tfrac{1}{2} [\Omega, \delta \Omega ] - \tfrac{1}{2} \left( K\tr \! (N - \tfrac{1}{2} K \delta \Omega) + (N - \tfrac{1}{2} K \delta \Omega)\tr \! K \right) = W,
\]
\[
   \delta \Omega + \tfrac{1}{2} \Omega \delta \Omega + \tfrac{1}{2} \delta \Omega \Omega - \tfrac{1}{2} \left( K\tr \! N - \tfrac{1}{2} K\tr \! K \delta \Omega + N\tr \! K - \tfrac{1}{2} \delta \Omega\tr K\tr \! K \right) = W,
\]
\[
   \delta \Omega + \tfrac{1}{2} \Omega \delta \Omega + \tfrac{1}{2} \delta \Omega \Omega - \tfrac{1}{2} K\tr \! N + \tfrac{1}{4} K\tr \! K \delta \Omega - \tfrac{1}{2} N\tr \! K + \tfrac{1}{4} \delta \Omega\tr K\tr \! K = W.
\]
We use the skew-symmetry $ \delta \Omega\tr = -\delta \Omega $ to get rid of the transpose,
\[
   \delta \Omega + \tfrac{1}{2} \Omega \delta \Omega + \tfrac{1}{2} \delta \Omega \Omega  + \tfrac{1}{4} K\tr \! K \delta \Omega  - \tfrac{1}{4} \delta \Omega K\tr \! K = W + \tfrac{1}{2} K\tr \! N + \tfrac{1}{2} N\tr \! K,
\]
collecting $ \delta \Omega $, we obtain
\[
   \delta \Omega + \tfrac{1}{2} \delta \Omega \Omega  - \tfrac{1}{4} \delta \Omega K\tr \! K + \tfrac{1}{2} \Omega \delta \Omega + \tfrac{1}{4} K\tr \! K \delta \Omega = W + \tfrac{1}{2} K\tr \! N + \tfrac{1}{2} N\tr \! K,
\]
\[
   \left( I_{p} +  \tfrac{1}{2} \Omega  + \tfrac{1}{4} K\tr \! K \right) \delta \Omega + \delta \Omega \left( \tfrac{1}{2} \Omega - \tfrac{1}{4}  K\tr \! K \right) = W + \tfrac{1}{2} K\tr \! N + \tfrac{1}{2} N\tr \! K,
\]
This is a Sylvester equation that can be solved with MATLAB's command \texttt{lyap} to obtain $\delta \Omega$. Then $\delta K$ can be found using \eqref{eq:delta_K}. Let the current iteration be indexed by $k$; then the tangent vector is updated as
\[
   \xi^{(k+1)} = \xi^{(k)} + Q \cdot
   \begin{bmatrix}
      \delta \Omega^{(k)} \\
      \delta K^{(k)}   
   \end{bmatrix}.
\]

%=========================================
% REFERENCES:
%=========================================
\bibliographystyle{aomalpha}

\begin{footnotesize}
   \bibliography{SSAF_preprint.bib}
\end{footnotesize}
%=========================================

\end{document}